\documentclass[11pt]{amsart}

\newtheorem{theorem}{Theorem}[section]

\theoremstyle{definition}

\theoremstyle{remark}

\numberwithin{equation}{section}

\usepackage{a4,a4wide}
\usepackage{leqno}
\usepackage{graphics}
\usepackage{amssymb}
\usepackage{amsmath}
\usepackage{amsthm}
\usepackage[latin1]{inputenc}
\usepackage{euscript}
\usepackage{palatino}
\usepackage{color}
\newenvironment{theo*}[2]{\smallskip\newline \noindent\textbf{Theorem}\;(#1)\;\newline \textit{#2}\medskip\\}

\newenvironment{th*}[1]{\smallskip\newline \noindent \textbf{Theorem}\;--\;\textit{#1}\smallskip}

\newcommand{\beq}{\begin{equation}}
\newcommand{\eeq}[1]{\label{#1}\end{equation}}

\def\O{{\Omega}}

\def\eps{{\epsilon}}

\def\k{{\mathcal{K}}}

\def\L{{\mathcal{L}}}

\def\M{{\mathcal{M}}}
\def\D{{\mathcal{D}}}

\def\R{{\mathbb{R}}}

\newcommand{\opk}[1]{\k[{#1}]}
\newcommand{\lb}[1]{\L_{_{#1}}}
\newcommand{\oplb}[2]{\L_{_{#2}}[{#1}]}
\newcommand{\opd}[1]{\D[{#1}]}

\newcommand{\opdb}[2]{\D_{_{#2}}[{#1}]}

\newcommand{\fdem}{\vskip 0.2 pt \qquad \qquad \qquad \qquad \qquad \qquad \qquad \qquad \qquad \qquad \qquad \qquad \qquad \qquad \qquad \qquad  $\square$  }

\def\tilde{\widetilde}

\begin{document}

\title{Singular measure as principal eigenfunction of some nonlocal operators }

\author{J{\'e}r{\^o}me Coville\\
~\\
\textit{\tiny  UR 546 Biostatistique et Processus Spatiaux\\
 INRA, Domaine St Paul Site Agroparc\\
  F-84000 Avignon\\
   France}}
\address{UR 546 Biostatistique et Processus Spatiaux\\
 INRA, Domaine St Paul Site Agroparc\\
  F-84000 Avignon\\
 France\\}
%\curraddr{INRA Avignon\\
%Unit\'e BioSP\\
%Site St Paul, Agroparc\\
%\\
%Avignon - France}
\email{jerome.coville@avignon.inra.fr}
\thanks{The author is supported by INRA Avignon   }
%\author{Lucas Rossi}
%\address{EHESS,}
%\email{}
%\thanks{Support information for the first author.}
%\author{Author Two}
%\address{Mathematical Research Section, School of Mathematical Sciences,
%Australian National University, Canberra ACT 2601, Australia}

%    General info

\date{February 4, 2013}



\begin{abstract}
In this paper, we are interested in the spectral properties of the generalised principal eigenvalue of some nonlocal operator. That is, we look for the  existence of some particular solution $(\lambda,\phi)$  of a nonlocal operator. 
$$\int_{\O}K(x,y)\phi(y)\, dy +a(x)\phi(x) =-\lambda \phi(x),$$
where $\O\subset\R^n$ is an open bounded connected set, $K$  a nonnegative kernel and $a$ is continuous. We prove that  for the generalised principal eigenvalue $\lambda_p:=\sup \{\lambda \in \R \, |\, \exists \, \phi \in C(\O), \phi > 0 \;\text{ so that }\; \oplb{\phi}{\O}+ a(x)\phi + \lambda\phi\le 0\}$ there exists always a solution $(\mu, \lambda_p)$ of the problem in the space of signed measure. Moreover $\mu$ a positive measure.  When $\mu$ is absolutely  continuous with respect to the Lebesgue measure, $\mu =\phi_p(x)$ is called the principal eigenfunction associated to $\lambda_p$.  In some simple cases, we exhibit some explicit singular measures that are solutions of the spectral problem.  
\end{abstract}

\maketitle
\section{\bf Introduction and Main results}
In this note we are interested in the spectral properties of some nonlocal operators. That is, we look for solution $(\phi,\lambda)$ of  
\begin{equation}
\int_{\O}K(x,y)\phi(y)dy   +(a(x)+\lambda) \phi(x)=0.
\end{equation}
%- \phi(x)\int_{\O}K(y,x)\,dy
where  $\O, K,$ and $a$ satisfy  the following assumptions : 
\begin{align*}
&\O\subset \R^n \quad \text{is a bounded open set}&\qquad(H1) \\ 
&K \in C(\R^n\times \R^n),\, K\ge 0,\,   \exists c_0 > 0,  \eps_0 > 0 \quad \text{such that} 
  \inf_{ x\in \O}\left( \inf_ {y\in B(x,\eps_0 )} K(x,y)\right) > c_0,&\qquad(H2)\\
&a \in C(\bar\O)\cap L^{\infty}(\O) &\qquad(H3)
\end{align*}
 
A typical  kernel which satisfies such assumptions is given by   
\begin{equation}
k(x,y) = J\left(\frac{x-y}{g(y)}\right)\frac{1}{g^n(y)},\label{pev-eq-ccem}
\end{equation}
where $J$ is a continuous positive probability density and the function $g$  is  bounded and positive.
Such kernel was introduced by Cortazar \textit{et al.} \cite{Cortazar2007a} in order to model a  non homogeneous dispersal process.

In the past few years much attention has been drawn to the study of nonlocal reaction-diffusion equations where such type of nonlocal operators is used  to model some long range effects. In this context,   the nonlocal operator  takes often  the form
\begin{equation}
\opdb{u}{\O}:=\int_{\O}K(x,y)u(y)\,dy - u(x)\int_{\O}K(y,x)u(y)\,dy, \label{gpev-eq-gene}
\end{equation}
  where $\O\subset \R^n$, $k \ge 0$ satisfies  $\int_{\R^n}
K(y,x)dy<\infty$ for all $x\in \R^n$; see among other references \cite{Bates2007,Cortazar2007a,Coville2010,Coville2012,Garcia-Melian2009,HMMV,Lutscher2005,SSN, Shen2012}.
In ecology these type of diffusion process has been widely used to describe the dispersal of a  population through its environment in the
following sense. As stated in  \cite{Fife1979,HMMV} if $u(y,t)$ is
thought of as a density at  a location $y$ at a time $t$ and $k(x,y)$
as the probability distribution of jumping from  a location $y$ to a
location $x$, then the rate at which the individuals from all other
places are arriving to the location $x$ is
$\int_{\O} k(x,y)u(y,t)\,dy.$
On the other hand, the rate at which the individuals are leaving the location $x$  is $-u(x,t)\int_{\O}k(y,x)u(y)\,dy $.
This formulation of the dispersal of individuals finds its justification in many ecological problems of seed dispersion; see for example \cite{CMS,Cl,DK,Medlock2003,SSN}.
%In such context, it is c see for example \cite{} where it is used  to model the dispersal of  the seeds   

The spectral properties  and in particular the existence of  principal eigenvalue have recently been investigated and some criteria for the existence of principal eigenvalue have been derived \cite{Coville2010,Coville2012,Kao2010,Shen2012}.  Namely,   the principal eigenvalue  $\lambda_p$ for the operator $\oplb{u}{\O}+a(x) u$      can be defined by the formula
 $$\lambda_p (\lb{\O} + a(x)) := \sup \{\lambda \in \R \, |\, \exists \, \phi \in C(\O), \phi > 0 \;\text{ so that }\; \oplb{\phi}{\O}+ a(x)\phi + \lambda\phi\le 0\}$$ and  as noted in \cite{Coville2010},  the condition $\frac{1}{\sup_{\bar\O}(a(x))-a(x)}\not \in L^1_{loc}(\bar \O)$ is sufficient to guarantees the existence of a continuous principal eigenfunction. Another useful criteria is 
\begin{theorem}[  \cite{Coville2012} ] \label{gpev-th-mup}
 There exists a positive continuous eigenfunction associated to 
$\lambda_p$ \textbf{if and only if} $\lambda_p(\lb{\O}+a)< -\sup_{\O}a$. 
\end{theorem}

On some examples, it has also been  observed  in \cite{Coville2010, Kao2010} that  $\lambda_p$ is not always an eigenvalue  in the space  $L^1(\bar\O)$.  %the condition $\frac{1}{\sup_{\bar\O}(a(x))-a(x)}\not \in L^1_{loc}(\bar \O)$ is sufficient to guarantees the existence of a continuous principal eigenfunction.

In this note, we study the properties of the principal eigenvalue $\lambda_p$ when no positive continuous eigenfunction exists i.e. when  $\lambda_p=-\sup_{\O}a$ and $\frac{1}{\sup_{\bar\O}(a(x))-a(x)} \in L^1(\bar \O)$.  In this situation, we construct  a positive measure,  solution of the equation
\begin{equation}
\int_{\O}K(x,y)d\mu(y)   +(a(x)+\lambda_p(\lb{\O}+a)) \mu=0.\label{gpev-eq1}
\end{equation}
%- \phi(x)\int_{\O}K(y,x)\,dy

This positive measure is a natural extension of the notion of principal eigenfunction associated to the generalised eigenvalue $\lambda_p$.
%solution, in the sense of distribution, of  the generalised principal eigenvalue problem
 %That is, we look for solution $(\mu,\lambda)\in \mb{\O}^+\times \R$ of  
%\begin{equation}
%\int_{\O}K(x,y)d\mu(y)   +(a(x)+\lambda) \mu=0.
%\end{equation}
%- \phi(x)\int_{\O}K(y,x)\,dy

%\begin{equation}
%\int_{\O}K(x,y)d\mu(y)   +(a(x)+\lambda_p(\lb{\O}+a)) \mu=0.. \label{gpev-eq1}
%\end{equation} 
%where $$\oplb{u}{\O}:= \int_{\O}K(x,y)u(y)dy   +a(x)u(x)$$
More precisely,
\begin{theorem}
Assume that  $K,a$ satisfies the assumptions $(H1-H3)$  and  so that $\frac{1}{\sup_{\bar \O}(a(x))-a(x)}\in L^1_{loc}(\bar \O)$.     Assume further that    $\lambda_p(\lb{\O}+a)=-\sup_{\O} a$ then there exists a  positive measure  $\mu \in \M^+(\O)$, solution to the equation \eqref{gpev-eq1}. Moreover,  if $\mu$ is singular with respect to the Lebesgue measure and $\#\O_0>1$ where $\O_0:=\{y\in \bar \O\, |\, a(y)=\sup_{\O}a\}$   then there exists infinitely many positive measures, solution of \eqref{gpev-eq1}.  
\end{theorem}

These measures are of great importance, since there are at the core of the analysis of propagating phenomena in heterogeneous nonlocal reaction diffusion like
\begin{equation}\label{gpev-eq-rd}
\frac{\partial u}{\partial t}= \opd{u}(x,t) +  u(x,t)(a(x)-u(x,t)).  
\end{equation}
The existence of a propagating speed seems to be  conditioned to the existence of continuous principal eigenfunction \cite{Coville2012,Shen2012}. 
They also  appear  naturally in the study of demo-genetic models such as 
\begin{equation}
\frac{\partial u}{\partial t}=u\left(a(x)-\int_{\O}a(y)u(y,t)\,dy\right) +\int_{\O}m(x,y)[u(y,t) -u(x,t)]\,dy, \label{gpev-eq-demogenet}
\end{equation}
where it is known that the solution $u(t,x)$ can blow up and converges to a solution of \eqref{gpev-eq1}, see for example \cite{Barles2008,Coville2013a,Jabin2011}.

%\int_{\O}J\left(\frac{x-y}{g(y)}\right)\frac{u(y)}{g^n(y)}\,dy -b(x)u=-\lambda u \quad \text{in} \quad \O.
%\end{equation}

%$|\O_0|=0$ where $\O_0:=\{y\in \bar \O\, |\, a(y)=\sup_{\O}a\}$  and

%Along this paper, with no further specifications, we will always make the following assumptions on   $\O$, $K$,  and $a$ : 
%\begin{align*}
%&\O\subset \R^n \quad \text{is a bounded open set}&\qquad(H1) \\ 
%&K \in C(\R^n\times \R^n),\, K\ge 0,\,   \exists c_0 > 0,  \eps_0 > 0 \quad \text{such that} 
%  \inf_{ x\in \O}\left( \inf_ {y\in B(x,\eps_0 )} K(x,y)\right) > c_0,&\qquad(H2)\\
%&a \in C(\bar\O)\cap L^{\infty}(\O) &\qquad(H3)
%\end{align*}
%where $C_c(\R^n)$ denotes the set of continuous functions with compact support.
%
%A typical example kernel satisfying such assumptions is   
%\begin{equation}
%k(x,y) = J\left(\frac{x-y}{g(y)}\right)\frac{1}{g^n(y)},\label{pev-eq-ccem}
%\end{equation}
%where $J$ is a continuous positive probability density and the function $g$  is  bounded and positive.
%This particular  type of diffusion kernel was introduced by Cortazar \textit{et al.} \cite{Cortazar2007a} in order to model a  non homogeneous dispersal process.  

\section{\bf Construction of the solution}
Let $x_0\in \bar \O$ be a point where the function $a$ achieves its maximum. Then \eqref{gpev-eq1} rewrites
\begin{equation}
\int_{\O}K(x,y)\phi(y)\,dy +(a(x)-a(x_0))\phi =0. \label{gpev-eq2}
\end{equation}
%As noted in \cite{Coville2010},  the condition $\frac{1}{a(x_0)-a(x)}\not \in L^1_{loc}(\bar \O)$ is sufficient to guarantees the existence of a continuous principal eigenfunction, so to construct a solution of \label{gpev-eq2} we have to assume that    
%$\frac{1}{a(x_0)-a(x)}\in L^1_{loc}(\bar \O)$. 

 Let us now introduce the  operator:
$$\opk{u}:=\int_{\O}K(x,y)\frac{u(y)}{a(x_0)-a(y)}\, dy.$$
Since $\frac{1}{a(x_0)-a(y)} \in L^1$ and $K \in C(\R^n\times\R^n)$ is non negative, the positive operator $\k$ is well defined in $C(\bar \O)$.  Moreover the operator $\k$ is  compact since $|\O_0|=0$ and $K$ is uniformly continuous in $\bar \O$.  Therefore its spectrum consists only on  isolated eigenvalue \cite{Brezis2010}.

Assume $1 \in \sigma(\k)$ then there exists  a positive $\psi \in L^1(\O)$ solution of \eqref{gpev-eq2}.
Indeed,   let $\lambda_1$ be the $\max\{|\lambda|\,|\,\lambda\in \sigma(\k)\}$ then  by the Krein-Rutman Theory there exists $\phi_1 \in C(\bar \O)$, $\phi_1>0$ so that  
 $$\opk{\phi_1}=\lambda_1\phi_1.$$  
Therefore there exists $\psi:=\frac{\phi}{(a(x_0)-a(x))} \in L^1$ so that for all $x\in \bar \O\setminus\O_0$
$$\int_{\O}K(x,y)\psi(y)\,dy+(a(x)-a(x_0))\psi=(\lambda_1 -1)\phi.$$ 
If $\lambda_1>1$, there exists $\eps>0$ and $\psi \ge 0$ so that for all  $x\in \bar \O\setminus\O_0$
$$\int_{\O}K(x,y)\psi(y)\,dy+(a(x)-a(x_0)-\eps)\psi=(\lambda_1 -1-\eps)\phi>0. $$
Thus all $x\in \bar \O\setminus\O_0$
$$\frac{1}{\eps+a(x_0)-a(x)}\int_{\O}K(x,y)\psi(y)\,dy> \psi(x)$$
and we get the contradiction
$$+\infty>\frac{\|K\|_{\infty}}{\eps}\int_{\O}\psi(y)\,dy \ge \limsup_{x \to x' \in \partial\O_0}\psi(x)=+\infty. $$
So $\lambda_1\le 1$ and by construction $\lambda_1=1$. 
Hence, in this situation there exists a positive measure $\phi_1$,  solution to \eqref{gpev-eq2}.

Now assume that $1\not \in  \sigma(\k)$. We will construct a positive measure $\mu$, solution of   \eqref{gpev-eq2}.
We  construct a solution of the form $\alpha \delta_{x_0} +f$ for some $\alpha \in \R^{+}$ and $f \in L^1(\O)$.  By introducing this Ansatz in \eqref{gpev-eq2} and after  a  straightforward computation we see that $f$ should satisfy 

\begin{equation}
\alpha K(x,x_0) + \int_{\O}K(x,y)f(y)\,dy  +(a(x)-a(x_0))f=0. \label{gpev-eq3}
\end{equation}

%We are then reduces to find $f$ in the right space.
By denoting $g=(a(x_0)-a(x))f$, we can rewrite \eqref{gpev-eq3}  as follows 
\begin{equation}
 \opk{g}-g =-\alpha K(x,x_0). \label{gpev-eq4}
\end{equation}
By assumption $1 \not \in \sigma(\k)$, so by the  Fredholm alternative the operator $\k -id$ is invertible and for any $h \in C(\bar \O)$ there exists a unique $g \in C(\bar \O)$ so that 
$$
 \opk{g}-g =h. 
$$
Moreover, $g\ge 0$ if $h\le 0$. Indeed, since 
 $$\opk{\phi_1}-\phi_1=(\lambda_1-1)\phi_1<0,$$
we have
\begin{align*}
\opk{g}-g &=\int_{\O}\tilde K(x,y)\phi_1(y) \frac{g(y)}{\phi_1(y)}\, dy -\frac{g(x)}{\phi_1(x)}\phi_1(x)\\
&=\int_{\O}\tilde K(x,y)\phi_1(y) \left[\frac{g(y)}{\phi_1(y)}- \frac{g(x)}{\phi_1(x)}\right]\, dy +(\lambda_1-1)\phi_1(x)\frac{g(x)}{\phi_1(x)}
\end{align*}
By denoting $w:=\frac{g}{\phi_1}$ we have for $w$
$$ \int_{\O}\tilde K(x,y)\phi_1(y) \left[w(y)- w(x)\right]\, dy +(\lambda_1-1)\phi_1(x)w(x)=h\le 0.$$
Since $\lambda_1<1$ and $\phi_1>0$ we easily conclude that $w$ cannot achieve a non positive  minimum without being constant.  Thus $w>0$ which implies that   $g>0$. 
% by using the definition of the principal eigenvalue we see that there exists $\phi_0>0$, $\phi_0 \in C(\bar \O)$ so that 
%$$\opk{\phi_0} -\phi_0 <0. $$
%Thus, $$h=\opk{g}-\frac{g}{\phi_0}\phi_0   $$
As a consequence,  for all $\alpha \in \R^{+,*}$, there exists a unique positive function $g_{\alpha}\in C(\bar \O)$  solution of  \eqref{gpev-eq4}. 
Let us denote $g_1$ the function obtained for $\alpha =1$. By construction, we have for any $\alpha \in \R$, $g_\alpha=\alpha g_1$. Therefore, if $\alpha\delta_{x_0}+\frac{g_{\alpha}}{a(x_0)-a(x)}$ is a solution of \eqref{gpev-eq2} then $$\alpha\delta_{x_0}+\frac{g_{\alpha}}{a(x_0)-a(x)}=\alpha\left( \delta_{x_0}+\frac{g_{1}}{a(x_0)-a(x)}\right) $$ and the constructed solution is unique up to normalisation. Furthermore any element $\mu$ of the linear set engendered by $\delta_{x_0}+\frac{g_{1}}{a(x_0)-a(x)}$, i.e.    
 $$ \mu\in Lin(\delta_{x_0}+\frac{g_{1}}{a(x_0)-a(x)}):=\left\{\alpha\left(\delta_{x_0}+\frac{g_{1}}{a(x_0)-a(x)}\right) \, | \alpha \in \R \right\}$$ is a signed-measure, solution of \eqref{gpev-eq2}. 

Observe that the above construction holds for any singular measure $\mu$ with its support contained in $\O_0$. 
Indeed, by replacing the Dirac measure at $x_0$ by $\mu$ in the above construction, we see that  $\mu+g$ is a positive measure solution of \eqref{gpev-eq2} where $g$ is the positive solution of  $$\opk{g} -g =-\int_{\O}K(x,y)d\mu(y).$$
 
 In particular when $\#\O_0>1$ the set of singular measures with support in $\O_0$ is not reduced to $\delta_{x_0}$  and we can construct a solution for any  points $x_0$  where the maximum of $a$ its achieves.  Thus any element $\nu$ in the linear space engendered by   $\delta_{x_0} + g_{x_0,1} $ and $ \delta_{x_1}+g_{x_1,1}$, i.e. $$\nu \in\left\{\alpha\left(\delta_{x_0}+\frac{g_{x_0,1}}{a(x_0)-a(x)}\right) +\beta \left(\delta_{x_1}+\frac{g_{x_1,1}}{a(x_1)-a(x)}\right) \, | \alpha, \beta \in \R \right\}$$ is a  signed-measure, solution of \eqref{gpev-eq2}. 
%Observe that the construction holds for any points where the function $a$ achieves its maximum. Thus when $\#\O_0>1$ the set of singular measure with support on $\O_0$ is not reduced to $\delta_{x_0}$ re exists at least two point $x_0$ and $x_1$ where the maximum its achieves and the linear space engendered by   $\delta_{x_0} + g_{x_0,1} $ and $ \delta_{x_1}+g_{x_1,1}$, $$\left\{\alpha\left(\delta_{x_0}+\frac{g_{x_0,1}}{a(x_0)-a(x)}\right) +\beta \left(\delta_{x_1}+\frac{g_{x_1,1}}{a(x_1)-a(x)}\right) \, | \alpha, \beta \in \R \right\}$$, is a set of sign-measure solution of \eqref{gpev-eq2}. 

 \fdem

\section{\bf Two simple examples }
\subsection*{First Example:} Let $\O=B_1(0)$ be the unit ball of $\R^3$ and let us consider  the following  eigenvalue problem
\begin{equation}
\rho\int_{\O}u+(a(x)+\lambda)u=0 \quad\text{ in }\quad \O \label{gpev-eqex}
\end{equation}
with $a(x)=1-\|x\|_2^2$ and $\rho>0$.
We can easily check that for $\rho> \frac{1}{4\pi}$  there exists a solution $(\phi_p,\lambda_p)$ to \eqref{gpev-eqex} with $\lambda_p<-sup_{x\in\O}(a(x))$ and $\phi_p \in C(\bar\O)$. Indeed, the function $\frac{\rho}{-\lambda_p- a(x)}$ with $\lambda_p$ so that  $\int_\O \frac{\rho}{-\lambda_p- a(x)}\,dx=1$ is the desired  solution. When $\rho=\frac{1}{4\pi}$ then \eqref{gpev-eqex} has still a positive solution $(\lambda_p,\phi_p)$ with $\lambda_p=-\max a$ and $\phi_p=\frac{\rho}{-\lambda_p- a(x)}$, the positive eigenfunction being now in $L^1(\O)$ and singular at $x_0$.
For $\rho<\frac{1}{4\pi}$ we can check that  there is no $L^1$ eigenfunction solution of  \eqref{gpev-eqex}. 
However, with $\lambda_p=-\max_{\O}a$, using the above construction for $\alpha=\frac{1}{\rho}-\int_{\O}\frac{dx}{a(x_0)-a(x)}$ we can see that the positive  measure
$$\alpha \delta_0+\frac{1}{a(x_0)-a(x)} $$
  is a solution of \eqref{gpev-eqex} .  
  Since $a$ achieves a unique maximum, this is the only positive measure, solution of \eqref{gpev-eqex}.
  %\fdem
\subsection*{Second Example:}  
Let $B_1(0)$ be the unit ball of $\R^2$ and let $\O=B_1(0)\times [0,1]$ be the unit cylinder of $\R^3$ and let us consider  the following principal eigenvalue problem
\begin{equation}
\rho\int_{\O}u+(a(x)+\lambda)u=0 \quad\text{ in }\quad \O \label{gpev-eqex2}
\end{equation}
with $a(x)=1-\sqrt{x_1^2+x_2^2}$ and $\rho>0$.
The maximum of $a(x)$ is achieved at any point $x_0$ of the symmetrical axis of the cylinder (i.e. $\O_0:=\{x_0:=(0,0,x_3)|\text{ with }  x_3\in (0,1)\}$) and we can check that 
$$\int_{\O}\frac{dx}{a(x_0)-a(x)}=\int_0^1\left(\int_{B_1(0)}\frac{dx_1dx_2}{\sqrt{x_1^2+x_2^2}} \right)dx_3=2\pi.$$
As in the previous example, for $\rho> \frac{1}{2\pi}$  there exists a solution $(\phi_p,\lambda_p)$ to \eqref{gpev-eqex2} with $\lambda_p<-sup_{x\in\O}(a(x))$ and $\phi_p \in C(\bar\O)$. When $\rho=\frac{1}{2\pi}$ then \eqref{gpev-eqex2} has still a positive solution $(\lambda_p,\phi_p)$ with $\lambda_p=-\max a$ and $\phi_p=\frac{\rho}{-\lambda_p- a(x)}$, the positive eigenfunction being now in $L^1(\O)$ and singular at any $x_0 \in \O_0$.
For $\rho<\frac{1}{2\pi}$ then there is no $L^1$ eigenfunction solution of  \eqref{gpev-eqex2}. But  with $\lambda_p=-\sup_{\O}a$ and for any $x_0\in \O_0$ we can then construct a positive measure $$\alpha \delta_{x_0}+\frac{1}{a(x_0)-a(x)}, $$  solution of \eqref{gpev-eqex2}.
On the interval $(0,1)$, let us consider a singular continuous measure $\mu$ like the devil staircase or the Ferenc Riez measure on $(0,1)$. Those two positive  measures are singular measure with their support in $\O_0$.   So from our construction $\mu +g_{\mu} $ is also a positive measure solution to \eqref{gpev-eqex2}. In these two cases, the singular measure does not have an atomic part.
%\fdem

From these two examples, we can see how  the set of measure that are solution to the spectral equation \eqref{gpev-eq1} can be rich. It is therefore necessary  to develop criteria to discriminate the pertinent solution.
  
% Here, since their is no principal eigenvalue, we know that $\lambda_p=-sup_{\O}a(x)$ and we are let to consider the solution of the following problem
%\begin{equation}
%\int_{\O}u+(a(x)-a(x_0))u=0 \label{gpev.eq2}
%\end{equation}
%
%To construct a solution, we seek a solution into the form $\alpha \delta_{x_0}+ \beta f$ where $\delta_{x_0}$ is the Dirac mass at $x_0$, $\alpha \in \R^+$ and $f$ is in $L^1(\O)$. 
%
%Plugging this Ansatz in  \eqref{qpev.eq2}, then we obtain
%\begin{equation*}
%\alpha + \int_{\O}f+(a(x)-a(x_0))f=0. \label{gpev.eq3}
%\end{equation*}
%By taking $f=\frac{1}{a(x_0)-a(x)}$ we end up with 
%\begin{equation*}
%\alpha + \int_{\O}\frac{dx}{a(x_0)-a(x)}-1 =0. \label{gpev.eq4}
%\end{equation*}
%Hence we have 
%$$\alpha=1-\int_{\O}\frac{dx}{a(x_0)-a(x)}.$$
%since $\int_{\O}\frac{dx}{a(x_0)-a(x)}<1$ because otherwise there exists a continuous  principal eigenfunction.
%The case $=$ should be studied.

%Observe that this computation works for any point where $a$ achieves a maximum. Therefore, if $a$ achieves a maximums at two distinct points $x_1$ and $x_2$ then the linear space engendered by the set of measures $\beta \delta_{x_1}+\gamma \delta_{x_2}+ \frac{1}{a(x_0)-a(x)}$, with  $\alpha = \beta+ \gamma$,   defines the   solutions of the problem. 
%If there is a continuum of point where the solution take its maximum 
%
%$$ \sum_{i} \alpha_i \delta_{x_i} + \sum_{j} \beta_j \delta_{\gamma_j} + ... + \sum_{n}\delta_{\Sigma_n} + \frac{1}{a(x_0)-a(x)}$$	 
% where $\Sigma_n$ are sets of dimension $n-1$. 
\section*{}
\bibliographystyle{plain}
\bibliography{measure.bib}

\end{document}